\documentclass[10pt,twoside,english]{amsart2000}
\advance\oddsidemargin by -1.0cm
\advance\evensidemargin by -1.0cm
\advance\topmargin by -1.0cm 
\usepackage{amssymb}
\usepackage{babel}
\usepackage{amstext2000}
\usepackage{amscd2000}   
\usepackage{epsfig}  
\usepackage{rotating}

\theoremstyle{plain}
\newtheorem{cor}{Corollary}[section]
\newtheorem{lem}{Lemma}[section]
\newtheorem{thm}{Theorem}[section]
\newtheorem{prop}{Proposition}[section]
\theoremstyle{definition}
\newtheorem{exa}{Example}[section]
\newtheorem{NB}{Remark}[section]

\newcommand{\bdm}{\begin{displaymath}}
\newcommand{\edm}{\end{displaymath}}
\newcommand{\be}{\begin{equation}}
\newcommand{\ee}{\end{equation}}
\newcommand{\ba}[1]{\begin{array}{#1}}
\newcommand{\ea}{\end{array}}

\newcommand{\btab}{\begin{tabular}}
\newcommand{\etab}{\end{tabular}}

\newcommand{\R}{\ensuremath{\mathbb{R}}}

\newcommand{\Ric}{\ensuremath{\mathrm{Ric}}}

\begin{document}
%
\thispagestyle{empty}
%
\date{\today}
\title{A relation between the Ricci tensor and the spectrum of the 
Dirac operator}
%
%
%
\author{Klaus-Dieter Kirchberg}
\address{ 
{\normalfont\ttfamily kirchber@mathematik.hu-berlin.de}\newline
Institut f\"ur Reine Mathematik \newline
Humboldt-Universit\"at zu Berlin\newline
Sitz: Rudower Chaussee 25\newline
D-10099 Berlin\\
Germany}
\thanks{This work was supported by the SFB 288 "Differential geometry
and quantum physics" of the Deutsche Forschungsgemeinschaft.}
\begin{abstract}
Using Weitzenb\"ock techniques on any compact Riemannian spin manifold we 
derive a general inequality depending on a real parameter and joining the 
spectrum of the Dirac operator with terms depending on the Ricci tensor and 
its first covariant derivatives. The discussion of this inequality yields 
vanishing theorems for the kernel of the Dirac operator $D$ and new lower 
bounds for the spectrum of $D^2$ if the Ricci tensor satisfies certain 
conditions.
\end{abstract}
\maketitle
\tableofcontents
\pagestyle{headings}
%
%
%
\section{Introduction}\noindent
%
%
In 1980 Th. Friedrich proved that, on any compact Riemannian spin manifold
$M^n$ of scalar curvature $R$ with $R_{\min} := \min \{ R(x) | x \in M\} >0$,
every eigenvalue $\lambda$ of the Dirac operator $D$ satisfies the
inequality
\begin{equation} \label{gl01}
\lambda^2 \ \ge \ \frac{n R_{\min}}{4(n-1)}
\end{equation}

(see \cite{2}). In special geometric situations, there are better estimations
than (\ref{gl01}) (see \cite{8}, \cite{9}). For example, in case of a compact K\"ahler
manifold $M^{2m}$ of complex dimension $m$ with positive scalar curvature $R$, we
have the estimate
\begin{equation} \label{gl02}
\lambda^2 \ \ge \ \left\{ \begin{array}{l}
\frac{m+1}{4m} \, R_{\min} \ \mbox{if $m$ is odd}\\[1em]
\frac{m}{4(m-1)} \, R_{\min} \ \mbox{if $m$ is even}
		      \end{array} \right. \ .
\end{equation}

Recently it was shown in \cite{5} and \cite{6} that the estimate (\ref{gl01})
can also be improved for such manifolds $M^n$ whose curvature tensor or Weyl tensor, respectively, is divergence free. It is well known that the 
curvature tensor $K$
of $M^n$ is divergence free if and only if the covariant derivative
$\nabla  \Ric$ of the Ricci tensor $\Ric$ has the property
\begin{equation} \label{gl03}
(\nabla_X \Ric)(Y) \ = \ (\nabla_Y \Ric)(X) \, , 
\end{equation}

where $X$ and $Y$ are arbritrary vector fields. In this paper we generalize
the results of \cite{5} in the sense that we do not make use of the 
condition (\ref{gl03}) here.

%
\section{The basic Weitzenb\"ock formula and first applications}
%
Let $M^n$ be any Riemannian spin manifold of dimension $n$ with Riemannian metric
$g$ and spinor bundle $S$. Then the twistor operator
\begin{displaymath}
\mathcal{D} : \Gamma (S) \to \Gamma (TM^n \otimes S)
\end{displaymath}

is defined by $\mathcal{D}{\psi} := X^k \otimes \mathcal{D}_{X_k}\psi$ and
$\mathcal{D}_X \psi := \nabla_X \psi + \frac{1}{n} X \cdot D \psi$. Here 
$D:= X^k \cdot \nabla_{X_k}$ denotes the Dirac operator, $(X_1 , \ldots, X_n)$ is
any local frame of vector fields and $(X^1, \ldots , X^n)$ is the associated
coframe given by $X^j := g^{jk} \cdot X_k$, where the $g^{jk}$ denote the 
components of the inverse of the matrix $(g_{jk})$ with $g_{jk} := g(X_j , X_k)$.
The curvature tensor $K$ of $M^n$ is defined by
\begin{displaymath}
K(X,Y)(Z) \ := \ \nabla_X \nabla_Y Z - \nabla_Y \nabla_X Z - \nabla_{[X,Y]}
Z \ . 
\end{displaymath}

$K$ and the curvature tensor $C$ of $S$ are related by
\begin{equation} \label{gl04}
C(X,Y) \psi \ = \ \frac{1}{4} X^k \cdot K(X,Y)(X_k) \cdot \psi \ . 
\end{equation}

For any vector field $X$, let $A_X$ denote the endomorphism of $S$ defined by
\begin{displaymath}
A_X \psi \ := \ (\nabla_{X_j} C)(X^j,X) \psi \stackrel{(4)}{=} \frac{1}{4}
 X^k \cdot (\nabla_{X_j} K)(X^j,X)(X_k) \cdot \psi \ . 
\end{displaymath}

The Bianchi identity implies
\begin{displaymath}
g\Big((\nabla_{X_j} K)(X^j, X)(Y),Z\Big) \ = \ g\Big((\nabla_Z \Ric)(Y) -
(\nabla_Y \Ric)(Z), X\Big) \ . 
\end{displaymath}

Thus we obtain
\begin{equation} \label{gl05}
A_X \psi \ = \ \frac{1}{4} \Big((\nabla_{X_k} \Ric)(X) \cdot X^k -
X^k \cdot (\nabla_{X_k} \Ric)(X)\Big) \cdot \psi \ . 
\end{equation}

Using this a simple calculation yields the identity
\begin{equation} \label{gl06}
X^k \cdot A_{X_k} \psi \ = \ - \frac{1}{4} (dR) \cdot \psi \ ,
\end{equation}

where $R:= \mathrm{tr} (\Ric)$ is the scalar curvature. Now, we consider the
family of differential operators
\begin{displaymath}
Q^t : \Gamma (S) \to \Gamma (TM^n \otimes S) 
\end{displaymath}

depending on $t \in \R$ and being locally defined by 
$Q^t \psi := X^k \otimes Q^t_{X_k} \psi$ and
\begin{displaymath}
Q^t_X \psi := \mathcal{D}_X \psi + t \cdot \Big(\Ric - \frac{R}{n}\Big)(X)
\cdot D \psi + t \cdot \Big(A_X - \frac{1}{4n} X \cdot dR\Big) \psi \ .
\end{displaymath}

Our first aim is to compute the function $|Q^t \psi|^2 := \langle
Q^t_{X_k} \psi , Q^t_{X^k} \psi \rangle$ for any eigenspinor $\psi$ of
the Dirac operator $D$. For any spinor field $\psi$, we have the general
identities
\begin{equation} \label{gl07}
X^k \cdot \mathcal{D}_{X_k} \psi \ = \ 0 \ , 
\end{equation}
\begin{equation} \label{gl08}
X^k \cdot \Big(\Ric - \frac{R}{n}\Big) (X_k) \cdot \psi \ = \ 0 \ , 
\end{equation}
\begin{equation} \label{gl09}
X^k \cdot \Big( A_{X_k} - \frac{1}{4n} X_k \cdot dR \Big) \psi \ = \ 0 \ ,
\end{equation}

which imply that the image of $Q^t$ is contained in the kernel of the Clifford 
multiplication, i.e., it holds that
\begin{equation} \label{gl10}
X^k \cdot Q^t_{X_k} \psi \ = \ 0 \ .
\end{equation}

Moreover, we have the well known relation
\begin{equation} \label{gl11}
| \mathcal{D} \psi |^2 \ = \ | \nabla \psi |^2 - \frac{1}{n} | D \psi |^2 \ . 
\end{equation}

\begin{lem} \label{lem2-1}
Let $\lambda$ be any eigenvalue of $D$ and let $\psi$ be any corresponding 
eigenspinor $(D \psi = \lambda \psi)$. Then, for all $t \in \R$,
the equation
\begin{equation} \label{gl12}
\begin{array}{l}
| Q^t \psi|^2 = |\nabla \psi|^2 - \frac{\lambda^2}{n} |\psi|^2 - 2t \Big(Re 
\langle \Ric (X^k) \nabla_{X_k} D \psi , \psi \rangle +\\[1em]
+ \mathrm{Re} \langle A_{X^k} \nabla_{X_k} \psi , \psi \rangle - 
\lambda^2 \frac{R}{n} |\psi|^2 \Big)+t^2\Big(\lambda^2 |\Ric - \frac{R}{n}|^2
 | \psi|^2- \\[1em]
- \lambda \langle ( \Ric (X_k) A_{X^k} + A_{X^k} \Ric (X_k)) \psi , \psi 
\rangle + \langle A_{X_k} \psi , A_{X^k} \psi \rangle -\\[1em]
- \frac{1}{16n} |dR|^2 |\psi |^2\Big) \ . 
\end{array}
\end{equation}

is valid.
\end{lem}

\begin{proof}
Using the identities (\ref{gl06}) - (\ref{gl09}) we calculate
\begin{eqnarray*}
|Q^t \psi |^2 &=& \langle \mathcal{D}_{X_k} \psi + t(\Ric - \frac{R}{n})(X_k) 
\psi +t (A_{X_k} - \frac{1}{4n} X_k \cdot dR) \psi \ , \\
&& \mathcal{D}_{X^k} \psi +t (\Ric - \frac{R}{n})(X^k) \psi +t( A_{X^k} - 
\frac{1}{4n} X^k \cdot dR) \psi \rangle \\
&=& | \mathcal{D} \psi |^2 +2t\lambda \mathrm{Re} \langle \nabla_{X_k} \psi , 
(\Ric - \frac{R}{n})(X^k) \psi \rangle +\\
&& +2t \mathrm{Re} \langle \nabla_{X_k} \psi , (A_{X^k} - \frac{1}{4n} X^k
\cdot dR) \psi \rangle +t^2 \lambda^2 |\Ric - \frac{R}{n}|^2 |\psi|^2+\\
&& 2t^2 \lambda \mathrm{Re} \langle \Ric (X_k) \cdot \psi , (A_{X^k} - 
\frac{1}{4n} X^k \cdot dR) \psi \rangle +\\
&& +t^2 \langle A_{X_k} \psi , (A_{X^k} - \frac{1}{4n} X^k \cdot dR) \psi
\rangle =\\
&=& | \mathcal{D} \psi |^2 - 2t \Big( \mathrm{Re} \langle \Ric (X^k)
\nabla_{X^k} D \psi , \psi \rangle - \mathrm{Re} \langle  \nabla_{X_k} \psi,
A_{X^k} \psi \rangle - \\
&& - \lambda^2 \frac{R}{n} |\psi |^2 - \frac{t}{2n} \lambda \mathrm{Re}
\langle \psi , dR \cdot \psi \rangle \Big)+ t^2 \Big( \lambda^2 | \Ric - 
\frac{R}{n} |^2 | \psi |^2 - \\
&& - 2 \lambda \mathrm{Re} \langle \psi, \Ric (X_k) \cdot A_{X^k} \psi 
\rangle - \lambda \frac{R}{2n} \mathrm{Re} \langle \psi , dR \cdot \psi 
\rangle +\\
&& + \langle A_{X_k} \psi , A_{X^k} \psi \rangle - \frac{1}{16n} |dR|^2 |\psi 
|^2 \Big) =\\
&=& |  \mathcal{D} \psi |^2 - 2t \Big( \mathrm{Re} \langle Ric (X^k)
\nabla_{X_k} D \psi, \psi \rangle + \mathrm{Re} \langle \psi, A_{X^k} 
\nabla_{X_k} \psi \rangle - \\
&& - \lambda^2 \frac{R}{n} |\psi |^2 \Big) +t^2 \Big( \langle A_{X_n} \psi , 
A_{X^k} \psi \rangle - \frac{1}{16n} |dR|^2 |\psi|^2 - \\
&& \lambda \langle  \psi , (\Ric (X_k) A_{X^k} + A_{X^k} \Ric (X_k)) \psi
\rangle + \lambda^2 | \Ric - \frac{R}{n}|^2 | \psi |^2 \Big) . 
\end{eqnarray*}

Inserting (\ref{gl11}) in the result of this calculation we obtain 
(\ref{gl12}).
\end{proof}

By Lemma 1.4 in \cite{5} and (\ref{gl05}), for all spinor fields $\psi$,
we have the identity
\begin{equation} \label{gl13}
\begin{array}{l}
\mathrm{Re} ( \langle \Ric (X^k) \cdot \nabla_{X_k} D \psi , \psi
\rangle ) + \mathrm{Re} (\langle A_{X^k} \nabla_{X_k} \psi , \psi
\rangle )=\\[1em]
=| \nabla D \psi |^2 - |(D^2 - \frac{R}{4} ) \psi |^2 - \frac{R}{4}
| \nabla \psi |^2 + \frac{1}{4} | \Ric |^2 | \psi |^2 +\\[1em]
+ \langle \nabla_{\Ric (X^k)} \psi , \nabla_{X_k} \psi \rangle - 
\mathrm{div} (X_{\psi}) \ , 
\end{array}
\end{equation}

where $X_{\psi}$ is the vector field defined by
\begin{displaymath}
X_{\psi} \ := \ \mathrm{Re} \Big( \langle (D^2 - \frac{R}{4} ) \psi,
\nabla_{X^k} \psi \rangle + \langle \nabla_{X_j} D \psi + \frac{1}{2}
\Ric (X_j) \cdot \psi , X^k \cdot \nabla_{X^j} \psi \rangle \Big) X_k \ . 
\end{displaymath}

Inserting (\ref{gl13}) into (\ref{gl12}) we obtain

\begin{thm} \label{thm-2-7}
Let $M^n$ be any Riemannian spin manifold and let $\lambda$ be any eigenvalue
 of the Dirac operator $D$. Then, for all $t \in \R$, any corresponding
eigenspinor $\psi$ satisfies the equation
\begin{equation} \label{gl14}
\begin{array}{l}
| Q^t \psi |^2 = | \nabla \psi |^2 - \frac{\lambda^2}{n} | \psi |^2 - 2t
\Big[ \langle \nabla_{\Ric (X_k)} \psi , \nabla_{X^k} \psi \rangle +\\[1em]
+ ( \frac{1}{4} | \Ric |^2 - \lambda^2 \frac{R}{n} ) | \psi |^2 +
(\lambda^2 - \frac{R}{4}) ( | \nabla \psi |^2 - (\lambda^2 - \frac{R}{4} )
| \psi |^2 ) - \mathrm{div} (X_{\psi}) \Big]+ \\[1em]
+ t^2 \Big[ \lambda^2 |\Ric - \frac{R}{n} |^2 | \psi |^2 - \lambda \langle
( \Ric (X^k) \cdot A_{X_k} + A_{X_k} \cdot \Ric (X^k)) \psi , \psi \rangle - 
\\[1em]
-  \langle A_{X_k} \cdot A_{X^k} \psi , \psi \rangle - \frac{1}{16n}
|dR|^2 | \psi |^2 \Big] \ . 
\end{array}
\end{equation}
\end{thm}

Equation (\ref{gl14}) is the basic Weitzenb\"ock formula of our paper.
In case of a divergence free curvature tensor, this formula simplifies
since we have $A_X =0$ then and $R$ is constant (compare with Theorem
1.5 in \cite{5}).\\
In the following we suppose that $M^n$ is compact. For any continuous function
$f$, we use the notation $f_{\min} (f_{\max})$ for the minimum (maximum) of
$f$ on $M^n$. Furthermore, let $\kappa$ denote the minimum of all eigenvalues
 of the Ricci tensor $\Ric$ on $M^n$. Then, for any spinor field $\psi$, we
have the inequality
\begin{equation} \label{gl15}
\kappa | \nabla \psi |^2 \ \le \ \langle \nabla_{\Ric (X_k)} \psi,
\nabla_{X^k} \psi \rangle \ . 
\end{equation}

In the following we use the Schr\"odinger-Lichnerowicz formula
\begin{equation} \label{gl16}
\nabla^* \nabla \ = \ D^2 - \frac{R}{4} \ . 
\end{equation}

\begin{lem} \label{lem-2-3}
Let $\lambda$ be any eigenvalue of the Dirac operator $D$. Then, for any 
corresponding eigenspinor $\psi$, there are the inequalities
\begin{equation} \label{gl17}
\begin{array}{c}
\displaystyle 
- \frac{R_{\max} - R_{\min}}{4} (\lambda^2 - \frac{R_{\min}}{4} )
\cdot \int_{M^n} | \psi |^2 \le \\[1em]
\displaystyle
\int_{M^n} (\lambda^2 - \frac{R}{4})(| \nabla \psi |^2 - (\lambda^2 - 
\frac{R}{4}) | \psi |^2 ) \le  \\[1em]
\displaystyle
\frac{R_{\max} - R_{\min}}{4} (\lambda^2 - \frac{R_{\min}}{4} )
\cdot \int_{M^n} | \psi |^2 \ . 
\end{array}
\end{equation}
\end{lem}

\begin{proof}
Using the formula (\ref{gl16}) we find
\begin{displaymath}
\int_{M^n} ( | \nabla \psi |^2 - ( \lambda^2 - \frac{R}{4})| \psi |^2 )
\ = \ 0 \ . 
\end{displaymath}

Hence, it follows that\\

$(*) \quad \quad \quad \displaystyle
\int_{M^n} (\lambda^2 - \frac{R}{4})(| \nabla \psi |^2 - (\lambda^2 - 
\frac{R}{4})| \psi |^2) = - \frac{1}{4} \int_{M^n} R\Big( | \nabla \psi |^2 
- (\lambda^2 - \frac{R}{4})| \psi |^2 \Big). $\\

Further, we have
\begin{eqnarray*}
&& \int_{M^n} R(|\nabla  \psi |^2 - ( \lambda^2 - \frac{R}{4}) | \psi
|^2 ) \le \\
& \le & R_{\max} \int_{M^n} | \nabla \psi |^2 - \int_{M^n} R( \lambda^2 - 
\frac{R}{4} ) | \psi |^2 \stackrel{(16)}{=}\\
&=& R_{\max} \int_{M^n} (\lambda^2 - \frac{R}{4}) | \psi |^2 - 
\int_{M^n} R( \lambda^2 - \frac{R}{4}) | \psi |^2 =\\
&=& \int_{M^n} (R_{\max} - R) ( \lambda^2 - \frac{R}{4}) | \psi |^2 \le
\int_{M^n} (R_{\max} - R)(\lambda^2 - \frac{R_{\min}}{4} )| \psi |^2 
\stackrel{(1)}{\le}\\
&\le &  \int_{M^n} (R_{\max} - R_{\min})(\lambda^2 - \frac{R_{\min}}{4})
| \psi |^2 = (R_{\max} - R_{\min})(\lambda^2 - \frac{R_{\min}}{4}) \cdot
\int_{M^n} | \psi |^2 . 
\end{eqnarray*}

This yields
\begin{displaymath}
- \frac{1}{4} \int_{M^n} R(| \nabla \psi |^2 - (\lambda^2 - \frac{R}{4})
|\psi|^2 ) \ge - \frac{R_{\max} - R_{\min}}{4} (\lambda^2 - \frac{R_{\min}}{4}
) \cdot \int_{M^n} | \psi |^2 . 
\end{displaymath}

Inserting this into $(*)$ we obtain the first one of the inequalities
(\ref{gl17}). The second one can be proved analogously.
\end{proof}

The endomorphisms of the spinor bundle $E:= - A_{X_k} \cdot A_{X^k}$ and
$T:= \Ric (X_k) \cdot A_{X^k} $ \linebreak $+ A_{X^k} \cdot \Ric (X_k)$ occur
on the right-hand side of formula (\ref{gl14}) . For any vector
field $X$, the endomorphism $A_X$ is antiselfadjoint
\begin{equation} \label{gl18}
(A_X )^* \ = \ - A_X \ . 
\end{equation}

Hence, $E$ and $T$ are selfadjoint
\begin{equation} \label{gl19}
E^* \ =  \ E \quad , \quad T^* \ = \ T \ . 
\end{equation}

Moreover, we see that $E$ is nonnegative
\begin{equation} \label{gl20}
E \ \ge \ 0 \ . 
\end{equation}

Let $\varepsilon$ denote the maximum of all eigenvalues of $E$ on $S$ and 
let $\tau$ be the minimum of the eigenvalues of $T$. Then, for any $\psi 
\in \Gamma (S)$, there are the inequalities
\begin{equation} \label{gl21}
\langle E \psi , \psi \rangle \ \le \ \varepsilon | \psi |^2 \ , 
\end{equation}
\begin{equation} \label{gl22}
- \langle T \psi , \psi \rangle \ \le \ - \tau | \psi |^2 \ . 
\end{equation}

In the following let $\lambda \ge 0$ be any nonnegative eigenvalue of $D$. 
Then, integrating equation (\ref{gl14}) using the inequalities
(\ref{gl15}), (\ref{gl17}), (\ref{gl21}) and (\ref{gl22}), in case
$t \ge 0$ and $\kappa \le 0$, we obtain
\begin{eqnarray*}
0 &\le& \int_{M^n} |Q^t \psi |^2 \le (\frac{n-1}{n} \lambda^2 - 
\frac{R_{\min}}{4} ) \cdot \int_{M^n} | \psi |^2 -\\
&& - 2t \Big( \kappa (\lambda^2 - \frac{R_{\min}}{4} ) \cdot
\int_{M^n} | \psi |^2 + ( \frac{1}{4} | \Ric |^2_{\min} - \lambda^2
\frac{R_{\max}}{n} ) \cdot \int_{M^n} | \psi |^2-\\
&& - \frac{R_{\max} - R_{\min}}{4} (\lambda^2 - \frac{R_{\min}}{4})
\cdot \int_{M^n} | \psi |^2 \Big) +\\
&& +t^2 \Big( \lambda^2 | \Ric - \frac{R}{n} |^2_{\max} - \lambda 
\tau + \varepsilon \Big) \cdot \int_{M^n} | \psi |^2 =\\
&=& \Big[ \frac{n-1}{n} \lambda^2 - \frac{R_{\min}}{4} - 2t \Big(
\frac{1}{4} | \Ric |^2_{\min} - \frac{R_{\max}}{n} \lambda^2 +\\
&& + (\kappa - \frac{R_{\max} - R_{\min}}{4}) (\lambda^2 - \frac{R_{\min}}{4}
) \Big) +\\
&& +t^2 (| \Ric - \frac{R}{n} |^2_{\max} \cdot \lambda^2 - \lambda \tau
+ \varepsilon ) \Big] \cdot \int_{M^n} | \psi |^2 \, . 
\end{eqnarray*}

Thus, in case $\kappa \le 0$, for all $t \ge 0$, it holds that
\begin{equation} \label{gl23}
\begin{array}{l}
\displaystyle
\lambda^2 - \frac{n R_{\min}}{4(n-1)} - 2t \Big( \frac{n}{4(n-1)} | \Ric 
|^2_{\min} - \frac{R_{\max}}{n-1} \lambda^2 + \\[1em]
\displaystyle
+ \frac{n}{n-1} ( \kappa - \frac{R_{\max} - R_{\min}}{4})(\lambda^2
- \frac{R_{\min}}{4} ) \Big)+ \\[1em]
\displaystyle
+ \frac{n}{n-1} t^2 \Big( | \Ric - \frac{R}{n} |^2_{\max} \cdot \lambda^2
- \tau \lambda + \varepsilon \Big) \ge 0 \, .
\end{array}
\end{equation}

In case $\kappa >0$, for all $t \ge 0$, it follows analogously that
\begin{equation}    \label{gl24}
\begin{array}{l}
\displaystyle
\lambda^2 - \frac{n R_{\min}}{4(n-1)} - 2t \Big( \frac{n}{4(n-1)} | 
\Ric|^2_{\min} - \frac{R_{\max}}{n-1} \lambda^2 +\\[1em]
\displaystyle
+ \frac{n \kappa}{n-1} (\lambda^2 - \frac{R_{\max}}{4})-
\frac{n(R_{\max}-R_{\min}}{4(n-1)} (\lambda^2 - \frac{R_{\min}}{4})
\Big)+\\[1em]
\displaystyle
+ \frac{n}{n-1} t^2 \Big( | \Ric - \frac{R}{n}|^2_{\max} \cdot \lambda^2
- \tau \cdot \lambda + \varepsilon \Big) \ge 0 \ . 
\end{array}
\end{equation}

From (\ref{gl01}) we know that $\ker (D) =0$ if the scalar curvature $R$ is
positive $(R_{\min} >0)$. Now, let us consider the case of $R_{\min} \le 0$.
This implies $\kappa \le 0$. Inserting $\lambda =0$ in (\ref{gl23}),
a simple discussion yields

\begin{thm} \label{thm-2-4}
Let $M^n$ be a compact Riemannian spin manifold with $R_{\min} \le 0$ such 
that the inequality
\begin{equation} \label{gl25}
| \Ric |^2_{\min} > R_{\min} \Big( \kappa - \frac{R_{\max}-R_{\min}}{4}
\Big) + 2 \sqrt{| R_{\min} | \cdot \varepsilon} 
\end{equation}

is satisfied. Then the kernel of the Dirac operator is trivial $(\ker
(D) =0)$.
\end{thm}

\begin{cor} \label{cor-2-5}
Let $M^n$ be a compact Riemannian spin manifold of constant scalar
curvature $R \le 0$ satisfying  the inequality
\begin{equation} \label{gl26}
| \Ric |^2_{\min} > R \kappa + 2 \sqrt{|R| \varepsilon} \ . 
\end{equation}

Then there are no harmonic spinors.
\end{cor}

\begin{cor} \label{cor-2-6}
Let $M^n$ be a compact Riemannian spin manifold satisfying the conditions
$R_{\min} =0$ and $| \Ric |_{\min} >0$. Then there are no harmonic spinors.
\end{cor}

We remark that Theorem \ref{thm-2-4} and Corollary \ref{cor-2-5} are
generalizations of Theorem \ref{thm-2-7} in \cite{5}. The case of
$R_{\min} >0$ gives rise to the question under which conditions the inequalities
(\ref{gl23}) and (\ref{gl24}), respectively, yield a better lower bound than
(\ref{gl01}). The answer can be given without determining the optimal
parameter $t$.\\

\underline{Case 1:} $R_{\min} >0 , \kappa \le 0$.\\
Inserting $\lambda := \sqrt{n R_{\min} /4(n-1)}$ into (\ref{gl23}) we
obtain a contradiction for some $t >0$ if and only if the coefficient
of $t$ is negative. But this is just the condition
\begin{equation} \label{gl27}
| \Ric |^2_{\min} > \frac{R_{\min}}{n-1} \Big( R_{\max} - \kappa
+ \frac{R_{\max} - R_{\min}}{4} \Big) \ .
\end{equation}

\underline{Case 2:} $R_{\min} >0 , \kappa >0$ . \\
The corresponding condition that can be derived from (\ref{gl24})
analogously is
\begin{equation} \label{gl28}
| \Ric |^2_{\min} > \frac{R_{\min}}{n-1} \Big( R_{\max} - \kappa +
\frac{R_{\max}-R_{\min}}{4} \Big) + \kappa (R_{\max} - R_{\min}).
\end{equation}

It is interesting to remark that in the inequalities (\ref{gl27}), 
(\ref{gl28})
the covariant derivatives of $\Ric$ do not appear. Moreover, in case of 
constant scalar curvature $R>0$, these conditions coincide and simplify
to
\begin{equation} \label{gl29}
| \Ric |^2_{\min} > \frac{R}{n-1} (R - \kappa) \ . 
\end{equation}

This generalizes the corresponding assertion in Section 2 of \cite{5}. \\[1em]

\section{The endomorphisms $E$ and $T$}

For further applications of the formulas (\ref{gl23}) and (\ref{gl24}),
it is convenient to write the endomorphisms
\begin{displaymath}
E = - A_{X_k} \cdot A_{X^k} \quad \mbox{and} \quad T= \Ric (X_k) \cdot
A_{X^k} + A_{X^k} \cdot \Ric (X_k)
\end{displaymath}

of the spinor bundle $S$ in a more suitable form.
\begin{prop} \label{prop-3-1}
There are the identities
\begin{equation} \label{gl30}
E= \frac{1}{4} | \nabla \Ric |^2 - \frac{1}{16} |dR|^2 + \frac{1}{8} [ 
\nabla_{X_j} \Ric , \nabla_{X_k} \Ric ] (X_l) \cdot X^j \cdot X^k \cdot
X^l \ , 
\end{equation}
\begin{equation} \label{gl31}
T = \frac{1}{2} [ \nabla_{X_k} \Ric , \Ric ] (X_l) \cdot X^k  \cdot X^l \ . 
\end{equation}
\end{prop}

\begin{proof}
We calculate
\begin{eqnarray*}
T &=& \frac{1}{4} \Big( \Ric (X^l) \cdot (\nabla_{X^k} \Ric)(X_l) \cdot
X^k - \Ric (X^l) \cdot X^k \cdot (\nabla_{X_k} \Ric )(X_l) \Big) +\\
&& \frac{1}{4} \Big( (\nabla_{X_k} \Ric )(X_l) \cdot X^k \cdot \Ric (X^l)-
X^k \cdot (\nabla_{X_k} \Ric)(X_l) \cdot \Ric (X^l) \Big) =\\
&=& \frac{1}{4} \Big( X^l \cdot (\nabla_{X_k} \Ric)(\Ric (X_l))\cdot
X^k - X^l \cdot X^k \cdot (\nabla_{X_k} \Ric )(\Ric (X_l)) +\\
&& + (\nabla_{X_k} \Ric)(\Ric (X_l)) \cdot X^k \cdot X^l - X^k \cdot
(\nabla_{X_k} \Ric)(\Ric (X_l)) \cdot X^l \Big) =\\
&=& \frac{1}{4} \Big( - (\nabla_{X_k} \Ric)(\Ric (X_l)) \cdot X^l \cdot X^k
- 2g( X^l, (\nabla_{X_k} \Ric \circ \Ric)(X_l)) X^k +\\
&& + X^l \cdot (\nabla_{X_k} \Ric)(\Ric (X_l)) \cdot X^k + 2g (X^k, (
\nabla_{X_k} \Ric \circ \Ric)(X_l)) X^l +\\
&& + (\nabla_{X_k} \Ric)(\Ric (X_l)) \cdot X^k \cdot X^l + (\nabla_{X_k}
\Ric)(\Ric (X_l)) \cdot X^k \cdot X^l +\\ 
&& + 2g ((\nabla_{X_k} \Ric \circ
\Ric)(X_l), X^k )X^l \Big)\\
&=& \frac{1}{2} \Big( (\nabla_{X_k} \Ric \circ \Ric)(X_l) \cdot X^k
\cdot X^l - (\nabla_{X_k} \Ric \circ \Ric)(X_l) \cdot X^l \cdot X^k-\\
&& -2 \mathrm{tr} (\nabla_{X_k} \Ric \circ \Ric )X^k + 2(\Ric \circ
\nabla_{X_k} \Ric)(X^k) \Big) =\\
&=& \frac{1}{2} \Big(( \nabla_{X_k} \Ric \circ \Ric)(X_l) \cdot X^k \cdot
X^l - (\nabla_{X_k} \Ric^2)(X_l) \cdot X^l \cdot X^k +\\
&& + (\Ric \circ \nabla_{X_k} \Ric )(X_l) \cdot X^l \cdot X^k -
\mathrm{tr} (\nabla_{X_k} \Ric^2)X^k +\\
&& +2 (\Ric \circ \nabla_{X_k} \Ric)(X^k) \Big) =\\
&=& \frac{1}{2} \Big( [ \nabla_{X_k} \Ric, \Ric ] (X_l) \cdot X^k \cdot
X^l - (\nabla_{X_k} \Ric^2)(X_l) \cdot X^l \cdot X^k - \\
&& - \mathrm{tr} (\nabla_{X_k} \Ric^2)X^k \Big) =\\
&=& \frac{1}{2} [ \nabla_{X_k} \Ric, \Ric ] (X_l) \cdot X^k \cdot X^l \ . 
\end{eqnarray*}

The latter equation is valid since $\nabla_X \Ric^2$ is selfadjoint
(symmetric) and, hence,
\begin{displaymath}
(\nabla_X \Ric^2)(X_l) \cdot X^l \ = \ - \mathrm{tr} (\nabla_X \Ric^2 ) \ . 
\end{displaymath}

Thus, we have (\ref{gl31}). Now, we prove (\ref{gl30}). It holds that
\begin{eqnarray*}
E &=& - \frac{1}{16} \Big( (\nabla_{X_j} \Ric)(X_l) \cdot X^j - X^j
(\nabla_{X_j} \Ric)(X_l)\Big) \Big( (\nabla_{X_k} \Ric)(X^l) \cdot X^k - 
X^k \cdot (\nabla_{X_k} \Ric)(X^l) \Big) =\\
&=& - \frac{1}{16} \Big( ( \nabla_{X_j} \Ric)(X_l) \cdot X^j \cdot 
(\nabla_{X_k} \Ric)(X^l)  \cdot X^k - (\nabla_{X_j} \Ric)(X_l) \cdot X^j
\cdot X^k \cdot (\nabla_{X_k} \Ric)(X^l) - \\
&& - X^j \cdot (\nabla_{X_j} \Ric)(X_l) \cdot (\nabla_{X_k} \Ric)(X^l) \cdot
X^k + X^j \cdot (\nabla_{X_j} \Ric)(X_l) \cdot X^k \cdot (\nabla_{X_k} \Ric)
(X^l) \Big) =\\
&=& - \frac{1}{16} \Big( ( \nabla_{X_j} \Ric \circ \nabla_{X_k} \Ric)(X_l)
\cdot X^j \cdot X^l \cdot X^k - ( \nabla_{X_j} \Ric \circ \nabla_{X_k} \Ric)
(X_l) \cdot X^j \cdot X^k \cdot X^l - \\
&& - X^j \cdot (\nabla_{X_j} \Ric \circ \nabla_{X_k} \Ric)(X_l) \cdot X^l 
\cdot X^k + X^j \cdot (\nabla_{X_j} \Ric \circ \nabla_{X_k} \Ric)(X_l) \cdot
X^k \cdot X^l \Big) =\\
&=& - \frac{1}{16} \Big( - 2 (\nabla_{X_j} \Ric \circ \nabla_{X_k} \Ric)
(X_l) \cdot X^j \cdot X^k \cdot X^l - 2(\nabla_{X_j} \Ric \circ
\nabla_{X_k} \Ric)(X^k) \cdot X^j\\
&& +2 X^j \cdot (\nabla_{X_j} \Ric \circ \nabla_{X_k} \Ric)(X_l) \cdot
X^k \cdot X^l + 2 X^j \cdot (\nabla_{X_j} \Ric \circ \nabla_{X_k} \Ric)(X^k)
\Big) =\\
&=& \frac{1}{8} \Big( (\nabla_{X_j} \Ric \circ \nabla_{X_k} \Ric)(X_l)
\cdot X^j \cdot X^k \cdot X^l + (\nabla_{X_j} \Ric \circ \nabla_{X_k} \Ric)
(X^k) \cdot X^j\\
&& + (\nabla_{X_j} \Ric \circ \nabla_{X_k} \Ric)(X_l) \cdot X^j \cdot X^k
\cdot X^l + 2g (X^j , (\nabla_{X_j} \Ric \circ \nabla_{X_k} \Ric)(X_l)
) X^k \cdot X^l -\\
&& - X^j \cdot (\nabla_{X_j} \Ric \circ \nabla_{X_k} \Ric)(X^k) \Big) =\\
&=& \frac{1}{8}\Big( 2(\nabla_{X_j} \Ric \circ \nabla_{X_k} \Ric)(X_l) \cdot
X^j \cdot X^k \cdot X^l +(\nabla_{X_j} \Ric \circ \nabla_{X_k} \Ric)(X^k) 
\Big) \cdot X^j+\\
&& +2g \Big( (\nabla_{X_k} \Ric \circ \nabla_{X_j} \Ric)(X^j), X_l) X^k \cdot
X^l - X^j \cdot (\nabla_{X_j} \Ric \circ \nabla_{X_k} \Ric)(X^k) \Big)=\\
&=& \frac{1}{8} \Big( - (\nabla_{X_k} \Ric \circ \nabla_{X^k} \Ric)(X_l) \cdot
X^l +[\nabla_{X_j} \Ric , \nabla_{X_k} \Ric] (X_l) \cdot X^j \cdot X^k \cdot
X^l+\\
&& +( \nabla_{X_j} \Ric \circ \nabla_{X_k} \Ric)(X^k)  \cdot X^j+ X^k
\cdot (\nabla_{X_k} \Ric \circ \nabla_{X_j} \Ric)(X^j) \Big)=\\
&=& \frac{1}{4} \mathrm{tr} (\nabla_{X_k} \Ric \circ \nabla_{X^k} \Ric)+
\frac{1}{8}[ \nabla_{X_j} \Ric , \nabla_{X_k} \Ric](X_l) \cdot X^j \cdot
X^k \cdot X^l-\\
&& - \frac{1}{4} g(X^j, (\nabla_{X_j} \Ric \circ \nabla_{X_k} \Ric)
(X^k)) =\\
&=& \frac{1}{4} | \nabla \Ric |^2 + \frac{1}{8} [ \nabla_{X_j} \Ric , 
\nabla_{X_k} \Ric ](X_l) \cdot X^j \cdot X^k \cdot X^l - \frac{1}{4} |
(\nabla_{X_k} \Ric)(X^k)|^2=\\
&=& \frac{1}{4} | \nabla \Ric |^2 - \frac{1}{16} |dR|^2 + \frac{1}{8}
[\nabla_{X_j} \Ric , \nabla_{X_k} \Ric ] (X_l) \cdot X^j \cdot X^k \cdot
X^l \ . 
\end{eqnarray*}
\end{proof}

Let $\Theta$ be the 3-form on $M^n$ defined by
\begin{eqnarray*}
\Theta (X,Y,Z) &:=& g ([\nabla_X \Ric , \Ric] (Y),Z)+g([\nabla_Z \Ric,
\Ric](X),Y)+\\
&& + g([\nabla_Y \Ric , \Ric](Z),X) \ . 
\end{eqnarray*}

\begin{cor} \label{cor-3-2}
The endomorphism $T$ acts on $S$ via Clifford multiplication by the 3-form 
$\Theta$, i.e., for all spinors $\psi$, we have the relation
\begin{equation} \label{gl32}
T \psi \ = \ \Theta \cdot \psi \ . 
\end{equation}
\end{cor}

\begin{proof}
By (\ref{gl31}), we have
\begin{displaymath}
T \psi = \frac{1}{2} g([\nabla_{X_k} \Ric, \Ric](X_l) , X_j) X^j 
\cdot X^k \cdot X^l \cdot \psi \ . 
\end{displaymath}

Thus, if $(X_1 , \ldots , X_n)$ is any local orthonormal frame, we obtain
\begin{eqnarray*}
T \psi &=& \frac{1}{2} \sum\limits_{\stackrel{j,k,l}{j=k}} g (
[\nabla_{X_k} \Ric , \Ric ] (X_l) , X_j) X_j \cdot X_k \cdot X_k \cdot 
\psi +\\
&& + \frac{1}{2} \sum\limits_{\stackrel{j,k,l}{k=l}} g([\nabla_{X_k}
\Ric, \Ric] (X_l), X_j) X_j \cdot X_k \cdot X_l \cdot \psi +\\
&& + \frac{1}{2} \sum\limits_{\stackrel{j,k,l}{k \not= j,l}}
g([\nabla_{X_k} \Ric, \Ric](X_l) , X_j) X_j \cdot
X_k \cdot X_l \cdot \psi =\\
&=& - \frac{1}{2} g([\nabla_{X_k} \Ric , \Ric](X_l), X^k)X^l \cdot
\psi - \frac{1}{2} g([\nabla_{X_k} \Ric, \Ric](X^k), X_j) X^j \cdot \psi
+\\
&& + \sum\limits_{\stackrel{j<l}{k \not= j,l}} g ([\nabla_{X_k} \Ric,
\Ric ](X_l), X_j) X_j \cdot X_k \cdot X_l \cdot \psi =\\
&=& \sum\limits_{\stackrel{j<l}{k \not= j,l}} g([\nabla_{X_k} \Ric,
\Ric](X_l) , X_j) X_j \cdot X_k \cdot X_l \cdot \psi =\\
&=& \sum\limits_{j<k<l} g([\nabla_{X_k} \Ric , \Ric](X_l) , X_j) X_j
\cdot X_k \cdot X_l \cdot \psi +\\
&& + \sum\limits_{k<j<l} g([\nabla_{X_k} \Ric , \Ric ] (X_l) , X_j) X_j
\cdot X_k \cdot X_l \cdot \psi +\\
&& + \sum\limits_{j<l<k} g([\nabla_{X_k} \Ric , \Ric](X_l) , X_j)
X_j \cdot X_k \cdot X_l \cdot \psi =\\
&=& \sum\limits_{j<k<l} (g([\nabla_{X_k}  \Ric , \Ric] (X_l) , X_j) -
g([\nabla_{X_j} \Ric , \Ric](X_l) , X_k)- \\
&& - g([\nabla_{X_l} \Ric , \Ric] (X_k) , X_j)) X_j \cdot X_k \cdot X_l
\cdot \psi =\\
&=& \sum\limits_{j<k<l} (g([\nabla_{X_k} \Ric , \Ric](X_l) , X_j )+
g ([\nabla_{X_j} \Ric , \Ric] (X_k), X_l)+\\
&& g([\nabla_{X_l} \Ric , \Ric] (X_j) , X_k)) X_j \cdot X_k \cdot X_l
\cdot \psi =\\
&=& \Theta \cdot \psi  \ . 
\end{eqnarray*}
\end{proof}

\section{Estimates for the first eigenvalue of the Dirac operator}

Let us introduce the notation
\begin{displaymath}
R_* := \left\{ \begin{array}{l}
R_{\min} \ \mbox{if $\kappa \le 0$}\\[1em]
R_{\max} \  \mbox{if $\kappa >0$}
	       \end{array} \ . 
\right.
\end{displaymath}

Then the inequalities (\ref{gl23}), (\ref{gl24}) can be written in the unified
form
\begin{equation} \label{gl33}
\alpha (t)  \lambda^2 - 2 \gamma (t)  \lambda \ge \beta (t) \quad \quad 
(\lambda \ge 0) \ , 
\end{equation}

where the functions $\alpha (t), \beta (t) , \gamma (t)$ 
depending on $t \ge 0$ are defined by
\begin{eqnarray*}
\alpha (t) &:=& 1+2t \frac{n}{n-1} \Big( \frac{R_{\max}}{n} - \kappa
+ \frac{R_{\max} - R_{\min}}{4} \Big) + \frac{n}{n-1} t^2 |\Ric - 
\frac{R}{n} |^2_{\max} , \\
\beta (t) &:=& \frac{n}{4(n-1)} \Big( R_{\min} +2t (| \Ric |^2_{\min}
- R_* \kappa + R_{\min} \frac{R_{max} - R_{\min}}{4} \Big) - 4t^2
\varepsilon \Big) , \\
\gamma (t) &:=& \frac{n \tau}{2(n-1)} t^2 \ . 
\end{eqnarray*}

By definition, we have $\alpha (t) \ge 1$ and the sign of $\gamma (t)$
depends on $\tau$. Moreover, we see that the inequality (\ref{gl33}) is
of interest only if the function $\beta (t)$ attains positive values for
some $t \ge 0$. Obviously, this is the case if $R_{\min} >0$. In case of
$R_{\min} \le 0$, this holds if the condition (\ref{gl25}) is
satisfied. Thus, from (\ref{gl33}) we immediately obtain 

\begin{thm} \label{thm-4-1}
Let $M^n$ be any compact $n$-dimensional Riemannian spin manifold and let
$\lambda \ge 0$ be any eigenvalue of the Dirac operator. Then, for all
$t \ge 0$ with $\beta (t) \ge 0$, there is the inequality
\begin{equation} \label{gl34}
\lambda \ge \frac{\sqrt{\alpha (t) \beta (t) +  \gamma (t)^2}+ \gamma (t)}
{\alpha (t)} = \frac{\beta (t)}{\sqrt{\alpha (t) \beta (t) + \gamma (t)^2}
- \gamma (t)} \ . 
\end{equation}
\end{thm}

\begin{thm} \label{thm-4-2}
Let $M^n$ be a compact Riemannian spin manifold with $\Theta =0$ and let
$\lambda$ be any eigenvalue of the Dirac operator. Then, for all $t \ge 0$,
we have the inequality
\begin{equation} \label{gl35}
\lambda^2 \ge \frac{\beta (t)}{\alpha (t)} \ . 
\end{equation}
\end{thm}

\begin{proof}
By Corollary \ref{cor-3-2}, $\Theta =0$ implies $T=0$. Thus, integrating
(\ref{gl14}) we obtain (\ref{gl35}). In this case, the sign of $\lambda$
plays no role.
\end{proof}

\begin{NB} \label{NB-1}
$M^n$ satisfies the condition $\Theta =0$ in the following special situations:
\begin{itemize}
\item[(i)] The covariant derivative of the Ricci tensor has the symmetry
property $(\nabla_X \Ric)(Y)=(\nabla_Y \Ric)(X)$, i.e., the curvature tensor
of $M^n$ is divergence free. This situation was investigated in \cite{5}.
\item[(ii)] The Ricci tensor commutes with its covariant derivatives, i.e.,
for all vector fields $X$, it holds that
\begin{equation} \label{gl36}
[  \nabla_X \Ric , \Ric ] =0 \ . 
\end{equation}

\noindent (\ref{gl36}) implies $\Theta =0$. This is an immediate consequence of
(\ref{gl31}) and (\ref{gl32}).
\item[(iii)] The Ricci tensor is recurrent, i.e., there is a 1-form $\eta$
on $M^n$ such that, for all vector fields $X$, the equation
\begin{equation} \label{gl37}
\nabla_X \Ric = \eta (X) \cdot \Ric
\end{equation}

\noindent is valid. Obviously, (\ref{gl37}) implies (\ref{gl36}). Thus, this
situation is a special case of\\situation (ii).
\end{itemize}
\end{NB}

\mbox{} 

\begin{NB} \label{NB-2}
In case the Ricci tensor has the property
\begin{equation} \label{gl38}
[ \nabla_X \Ric , \nabla_Y \Ric ] =0 \ , 
\end{equation}

i.e., if any two covariant derivatives of $\Ric$ commute, it follows
from (\ref{gl30}) that the endomorphism $E$ simplifies to the function
\begin{equation} \label{gl39}
E= \frac{1}{4} | \nabla \Ric |^2 - \frac{1}{16} |dR|^2 \ . 
\end{equation}

Hence, in this case we have
\begin{equation} \label{gl40}
\varepsilon = \Big( \frac{1}{4} | \nabla \Ric |^2 - \frac{1}{16} |dR|^2
\Big)_{\max} \ . 
\end{equation}

We see that (\ref{gl37}) implies (\ref{gl38}). Moreover, if $\Ric$ has
pairwise different eigenvalues on $M^n$, then (\ref{gl36}) implies
(\ref{gl38}).
\end{NB}

Calculating the maximum of the function $\beta (t) /  \alpha (t)$ by
Theorem \ref{thm-4-2} we obtain

\begin{cor} \label{cor-4-3}
Let $M^n$ be a compact Riemannian spin manifold satisfying the conditions
$R_{\min} =0 , | \Ric|_{\min} >0$ and $[\nabla_X \Ric , \Ric] =0$. Then,
for every eigenvalue $\lambda$ of the Dirac operator, we have the estimate
\begin{equation} \label{gl41}
\lambda^2 \ge \frac{n}{8(n-1)} \cdot \frac{| \Ric |^2_{\min}}{a+b+
\sqrt{a^2 +ab+c}} \ , 
\end{equation}

where the constants $a,b,c$ are defined by 
\begin{eqnarray*}
a&:=& \frac{\varepsilon}{| \Ric |^2_{\min}} \ , \\
b &:=& \frac{n}{2(n-1)} \Big( \frac{n+4}{4n} R_{\max} - \kappa \Big) \ , \\
c &:=& \frac{n}{4(n-1)} | \Ric - \frac{R}{n} |^2_{\max} \ . 
\end{eqnarray*}
\end{cor}

Now, we give some simple examples.

\begin{exa}
Let us consider the Riemannian product $M^4 (r, \rho):=S^2(r) \times
T^2 (\rho)$, where $S^2 (r) \subset \R^3$ is the standard sphere of radius
$r>0$ and $T^2 (\rho) \subset \R^3 (\rho >0)$ the standard torus defined
by $\mathfrak{x} :[0,2 \pi] \times [0,2 \pi] \to \R^3$,
$\mathfrak{x} (u,v):= (\rho (2+ \cos u)  \cos v, \rho (2+ \cos u) \sin v, 
\rho 
\sin u)$. Then the Ricci tensor $\Ric$ and the scalar curvature $R$ of
$M^4 (r,\rho)$ are given by\\

$\displaystyle (*) \hfill \Ric = \frac{1}{r^2} p_S + \frac{\cos u}{\rho^2
(2+ \cos u)} p_T \ , $ \hfill \mbox{}\\[1em]

$\displaystyle (2*) \hfill R= \frac{2}{r^2} + \frac{2 \cos u}{\rho^2 (2+
\cos u)} \ , $ \hfill \mbox{}\\

where $p_S$ and $p_T$ denote the projections on the tangent spaces of 
$S^2 (r)$ and $T^2 (\rho)$, respectively. This implies\\

$(3*) \hfill \displaystyle  | \Ric |^2 = \frac{2}{r^4} + \frac{2 \cos^2 u}
{\rho^4 (2+ \cos u)^2} \ , $  \hfill \mbox{} \\[1em]

$(4*) \hfill \displaystyle | \Ric - \frac{R}{4}|^2 = \Big( \frac{1}{r^2} - 
\frac{\cos u}{\rho^2(2+ \cos u)} \Big)^2$ , \hfill \mbox{}\\[1em]

$(5*) \displaystyle  \hfill \nabla_X \Ric = - \frac{2 \sin u}{\rho^2
(2+ \cos u)^2} du(X) p_T$ , \hfill \mbox{}\\[1em]

$(6*) \displaystyle \hfill \frac{1}{4} | \nabla \Ric |^2 - \frac{1}{16}
|dR|^2 = \frac{\sin^2 u}{\rho^6 (2+ \cos u)^4}$ \ . \hfill \mbox{}\\[1em]

Hence, we obtain\\

$(7*)  \displaystyle \hfill | \Ric |^2_{\min} = \frac{2}{r^4}$ \hfill
\mbox{}\\[1em]

$(8*) \displaystyle \hfill R_{\min} = 2 \Big( \frac{1}{r^2} - 
\frac{1}{\rho^2} \Big) , R_{\max} = 2 \Big( \frac{1}{r^2} + 
\frac{1}{3\rho^2} \Big)$ , \hfill \mbox{}\\[1em]

$(9*) \displaystyle \hfill | \Ric - \frac{R}{4} |^2_{\max} = \Big(
\frac{1}{r^2} + \frac{1}{\rho^2} \Big)^2 $ , \hfill \mbox{}\\[1em]

$(10*) \displaystyle \hfill \kappa = - \frac{1}{\rho^2} $ . \hfill
\mbox{} \\[1em]

From $(*)$ and $(5*)$ we see that there are the relations
\begin{displaymath}
[\nabla_X \Ric , \Ric ]=0 \quad , \quad [\nabla_X \Ric , \nabla_Y \Ric
]=0 \ . 
\end{displaymath}

Thus, by Proposition \ref{prop-3-1}, Remark \ref{NB-2} and $(6*)$,
we find here\\

$(11*) \displaystyle \hfill \tau =0 \quad , \quad \varepsilon =
\frac{3+2 \sqrt{3}}{36 \rho^6} $  . \hfill \mbox{}\\[1em]
\end{exa}
 
\underline{The case of $r= \rho (R_{\min} =0):$}\\

In this case, the suppositions of Corollary \ref{cor-4-3} are satisfied.
By $(7*)$ - $(11*)$, the constants $a,b,c$ are given by
\begin{displaymath}
a= \frac{3+ 2 \sqrt{3}}{72 r^2} \quad , \quad b= \frac{14}{9r^2} \quad
, \quad c= \frac{4}{3 r^4} \ . 
\end{displaymath}

Inserting this into (\ref{gl41}) we obtain the estimate
\begin{displaymath}
\lambda^2 \ge 0,116 \, r^{-2} \ . 
\end{displaymath}

\underline{The case of $r< \rho (R_{\min} >0)$:}\\

By $(7*)$ - $(10*)$, (\ref{gl27}) is equivalent to
\begin{displaymath}
\frac{1}{r^4} > \frac{1}{3 \rho^2} \Big( \frac{1}{r^2} - \frac{7}{\rho^2}
\Big) \ . 
\end{displaymath}

For $r < \rho$, this inequality is always satisfied. Hence, for all manifolds
$M^4 (r,\rho)$ with $r <\rho$, (\ref{gl35}) yields a better estimate than
(\ref{gl01}). For example, let us consider the special case that $\rho = \frac{r}
{3} \sqrt{10}$. Then, by (\ref{gl35}) and $(7*)$ - $(11*)$, we find\\

$(12*) \displaystyle \hfill \lambda^2 \ge 0,156 \, r^{-2} $ , \hfill
\mbox{}\\

whereas (\ref{gl01}) yields the estimate
\begin{displaymath}
\lambda^2 \ge \frac{1}{15} r^{-2} = 0,0 \bar{6} \, r^{-2} \ . 
\end{displaymath}

Since $M^4 (r, \rho)$ is K\"ahler, we can also apply the estimate
(\ref{gl02}). For $\rho = \frac{r}{3} \sqrt{10}$, (\ref{gl02}) and
$(8*)$ yield 
\begin{displaymath}
\lambda^2 \ge 0,1  \, r^{-2} \ . 
\end{displaymath}

Thus, the estimate $(12*)$ obtained by Theorem \ref{thm-4-2} is the better
one.\\

\underline{The case of $r>\rho (R_{\min} <0)$:}\\

For example, let us consider the special case that $\rho = \frac{3}{\sqrt{10}} r$,
i.e., $R_{\min} = - \frac{2}{9} r^{-2}$. Then one finds that the
condition (\ref{gl25}) is satisfied. Thus, Theorem \ref{thm-2-4} implies
that there are no harmonic spinors on $M^4 (r, \frac{3}{\sqrt{10}} r)$.
Calculating the maximum of the corresponding function $\beta (t) / \alpha 
(t)$, our Theorem \ref{thm-4-2} yields the estimate
\begin{displaymath}
\lambda^2 \ge 0,061 \, r^{-2} \ . 
\end{displaymath}

\begin{exa}
For $\rho >1$, we consider the Riemannian product $M^6 (\rho) := S^2
\times F^2 \times T^2 (\rho)$, where $S^2 \subset \R^3$ is the unit 
standard sphere, $F^2$ any compact Riemannian surface of constant Gaussian
curvature $-1$ and $T^2 (\rho) \subset  \R^3$ the torus defined in
Example 4.1. Then the Ricci tensor and scalar curvature of $M^6 
(\rho)$ are given by
\begin{displaymath}
\Ric = p_S - p_F + \frac{\cos u}{\rho^2 (2+ \cos u)} p_T \ , 
\end{displaymath}
\begin{displaymath}
R= \frac{2 \cos u}{\rho^2 (2+ \cos u)} \ , 
\end{displaymath}

where $p_S , p_F, p_T$ denote the projections on the tangent spaces of $S^2,
F^2$ and $T^2 (\rho)$, respectively. This implies\\

$(*) \displaystyle \hfill | \Ric |^2_{\min} =4 $  , \hfill \mbox{}\\[1em]

$(2*) \displaystyle \hfill R_{\min} = - \frac{2}{\rho^2} \quad , \quad
R_{\max} = \frac{2}{3 \rho^2}$ , \hfill \mbox{}\\[1em]

$(3*) \displaystyle \hfill | \Ric - \frac{R}{6} |^2_{\max} = 4 \Big( 1+
\frac{1}{3  \rho^4} \Big)$ . \hfill \mbox{}\\

Moreover, here we also have $[\nabla_X \Ric , \Ric ]=0$ and
$[\nabla_X \Ric , \nabla_Y \Ric ]=0$. Hence, as in Example 1 we obtain\\

$(4*) \displaystyle \hfill \tau =0 \quad , \quad \varepsilon = \frac{3+2 
\sqrt{3}}{36 \rho^6}$ . \hfill \mbox{}\\

Furthermore, $\rho >1$ implies\\

$(5*) \displaystyle \hfill \kappa = - 1$ . \hfill \mbox{}\\

By $(*)$ - $(5*)$, we find
\begin{eqnarray*}
\alpha (t) &=& 1+ (\frac{12}{5} + \frac{28}{15} \rho^{-2} )t +(
\frac{24}{5} + \frac{8}{5} \rho^{-4}) t^2 , \\
\beta (t) &=& - \frac{3}{5} \rho^{-2} +(\frac{12}{5} - \frac{6}{5} 
\rho^{-2} - \frac{4}{5} \rho^{-4} )t - \frac{3 +2 \sqrt{3}}{30} 
\rho^{-6} t^2
\end{eqnarray*}

and (\ref{gl25}) is equivalent to
\begin{displaymath}
\rho^4 - \frac{1}{2} \rho^2 > \frac{1}{3} + \frac{1}{12} \sqrt{6 + 4 
\sqrt{3}} \ . 
\end{displaymath}

The last inequality is valid for $\rho > \rho_0 \approx 1,04113$. 
Thus, by Theorem \ref{thm-2-4}, $M^6 (\rho)$ admits
no harmonic spinors for $\rho > \rho_0$. 
Moreover, we are interested in the maximum 
$f(\rho)(t_{\max})$ of the function $f(\rho)(t):= \beta (t) / \alpha
(t)$ for $t \ge 0$ subject to $\rho$ describing the lower bound of 
$\lambda^2$ on $M^6 (\rho)$. We obtain the following picture:

\mbox{} \hspace{2cm}
\begin{minipage}[b]{.46\linewidth}
\epsfig{figure=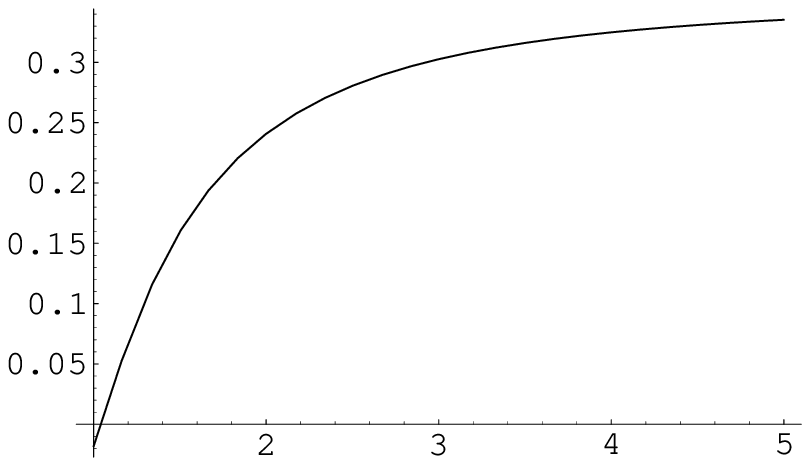}
\end{minipage}

The lower bound is a bounded monotone increasing function of $\rho$ \\with
$\lim\limits_{\rho \to \infty} f(\rho)(t_{\max}) \approx 0,354$. \\For
example, on $M^6 (2)(R_{\min} = - 1/2)$ we have the estimate
$\lambda^2 \ge 0,2407$. 
\end{exa}

\begin{NB} \label{NB-3}
The manifolds $M^4 (r, \rho)$ and $M^6 (\rho)$ of Example 1 and Example 2, 
respectively, do not satisfy the condition (\ref{gl03}). These two
examples are special cases of the following more general situation where
the geometric suppositions (\ref{gl25}), (\ref{gl27}) and (\ref{gl28}),
respectively, can be
checked easily: Let us consider a Riemannian product of the form 
$M^n := M^{n_1}_1 \times \ldots \times M^{n_p}_p$, where $M^{n_k}_k
(k=1 , \ldots , p)$ is a compact Riemannian spin manifold of dimension
$n_k$ such that $M^{n_k}_k$ is Einstein if $n_k \ge 3$. Then the Ricci tensor
$\Ric$ of $M^n$ satisfies the conditions (\ref{gl36}) and (\ref{gl38}), but
not the condition (\ref{gl03}) in general. Thus, here we obtain
\begin{displaymath}
\tau =0 \quad , \quad \varepsilon = \frac{1}{16} \cdot
\sum\limits^p_{k=1} |dR_k|^2_{\max} \ , 
\end{displaymath}

where $R_k$ denotes the scalar curvature of $M^{n_k}_k$ and, hence, $dR_k =0$ 
if $n_k \not= 2$. Moreover, we have
\begin{eqnarray*}
| \Ric |^2_{\min} &=& \sum\limits^p_{k=1} \frac{(R^2_k)_{\min}}{n_k} , \\
R_{\min / \max} &=& \sum\limits^p_{k=1} (R_k)_{\min / \max} , \\
\kappa &=& \min \Big\{ \frac{(R_k)_{\min}}{n_k} | k=1, \ldots , p \Big\} , \\
| \Ric - \frac{R}{n} |^2_{\max} &=& \Big( \sum\limits^p_{k=1} n_k (
\frac{R_k}{n_k} - \frac{R}{n} )^2 \Big)_{\max}
\end{eqnarray*}

with $R=R_1 + \ldots + R_p$ and $n=n_1 + \ldots + n_p$. Thus, we
see that, in this situation, all the essential numbers can be derived from the
scalar curvatures $R_1 , \ldots , R_p$ and the dimensions $n_1 , \ldots , 
n_p$.
\end{NB}

\begin{NB}
Up to now we have no answer to the question if there exist compact 
Riemannian spin manifolds that realizes the limiting case of the inequality
(\ref{gl34}) in case of an optimal parameter $t= t_0 >0$.
\end{NB}


%

\end{document}